\documentclass[twoside]{ima} 
\usepackage{amsfonts} 

\newcommand \ubar   	{{\overline u}}
\newcommand \vbar   	{{\overline v}} 
\newcommand \del		\partial  
\newcommand \auth {\textsc}
\newcommand \la 		\langle
\newcommand \ra 		\rangle
\newcommand \be   	{\begin{equation}}
\newcommand \ee   	{\end{equation}}
\newcommand \sgn   {\mbox{sgn }}
\newcommand \RR      {\mathbb{R}}
\newcommand \Hcal   {\mathcal{H}}
\newcommand \omegab {\overline \omega} 
\newcommand \Lcal {\mathcal L}  
\newcommand \dive {\textrm{div}_\omega\hskip.06cm}
\newcommand \Sphere {\mathbf S^2}
\newcommand \qbf {\mathbf q}

\begin{document} 
\markboth{P{\scriptsize HILIPPE} G. L{\scriptsize E}F{\scriptsize LOCH}}{HYPERBOLIC CONSERVATION LAWS ON SPACETIMES}          
\title{HYPERBOLIC CONSERVATION LAWS ON SPACETIMES} 
 \author{P{\scriptsize HILIPPE} G. L{\scriptsize E}F{\scriptsize LOCH}\thanks{
Laboratoire Jacques-Louis Lions \& Centre National de la Recherche Scientifique,
Universit\'e Pierre et Marie Curie (Paris 6), 4 Place Jussieu,  75252 Paris, France.
 Blog: {\tt philippelefloch.wordpress.com.}
\newline 
The author was supported by the Agence Nationale de la Recherche via the grant ANR~06-2-134423.   
{\tt To appear in:} Summer Program on ``Nonlinear Conservation laws and applications'', edited by A. Bressan, G.-Q. Chen, M. Lewicka, and D. Wang, July 13--31, 2009, IMA publication, 2010. 
}}   

\maketitle   
             
\begin{abstract}      
We present a generalization of Kruzkov's theory to manifolds.
Nonlinear hyperbolic conservation laws are posed on a differential $(n+1)$-manifold, 
called a spacetime, and the flux field is defined as a field of $n$-forms depending on a parameter. 
The entropy inequalities take a particularly simple form 
as the exterior derivative of a family of $n$-form fields. 
Under a global hyperbolicity condition on the spacetime, which 
allows arbitrary topology for the spacelike hypersurfaces of the foliation, 
we establish the existence and uniqueness of an entropy solution
to the initial value problem, and we derive a geometric version of the 
standard $L^1$ semi-group property. We also discuss an alternative framework in which the flux field 
consists of a parametrized family of vector fields.              
\end{abstract} 
 
\AMSMOS 
{Primary: 35L65. Secondary: 76L05, 76N.}  
\endAMSMOS

 
\section{Introduction} 
 
Hyperbolic conservation laws on manifolds arise in many applications to continuum physics.
Specifically, the shallow water equations for perfect fluids posed on the sphere with prescribed topography
provide an important model describing geophysical flows and exhibiting complex wave structures (Rosby waves, etc.) 
with several distinct scales.  
In \cite{ABL}--\cite{ 
BFL} and \cite{GL,LNO,LO1,LO2}, the author together with several collaborators recently initiated a research program 
and posed the foundations for the analysis and numerical approximation 
of hyperbolic conservation laws on manifolds. These equations form a simplified, yet challenging, mathematical model 
for studying the propagation of nonlinear waves, including shock waves, and their interplay with the background geometry.  


\section{Hyperbolic conservation laws based on vector fields} 
  
The standard shock wave theory posed on the Euclidian space covers nonlinear hyperbolic equations of the form 
$$
\del_t u +  \sum_{j=1}^n \del_j f^j(u) = 0, \qquad u=u(t,x) \in \RR, \, \, t \geq 0, \, x \in \RR^n, 
$$  
where $f^j: \RR \to \RR$ are given functions, 
and, more generally, applies to balance laws with variable coefficients and source-terms. 
Pioneering works by Lax, Oleinik, and others in the 50's and 60's set the foundations about discontinuous entropy 
solutions to such equations. Later, well-posedness theorems within the class of entropy solutions were established by 
Conway and Smoller (1966), Volpert (1967) (for solutions with bounded variation), 
Kruzkov (1970) (for bounded solutions), and finally DiPerna (1984) (for measure-valued solutions). 
Based on this mathematical theory, active research followed in the 80's and 90's 
that led to the development of robust and high-order accurate, shock-capturing schemes. 

The generalization of the above theory to {\sl manifolds} has only recently received particular attention, 
driven by problems arising in geophysical science and general relativity. In the present section, we 
consider 
a class of conservation laws based on parametrized vector fields,
 and we suppose that $M$ is a compact, oriented $n$-dimensional manifold,
endowed with
an $L^\infty$ volume form such that, in a finite atlas of local coordinate charts $(x^j)$ and for uniform constants 
$c \leq C$, 
$$
\omega = \omegab \, dx^1 \cdots dx^n, \qquad 0 < c \leq \omegab \leq C. 
$$ 
Recall that the divergence of an $L^\infty$ vector field $X$ is defined in the sense of distributions by 
$$
\la \dive X, \theta \ra :=  \int_M (d\theta) ( X) \, \omega
$$
for all test-function $\theta: M \to \RR$. 

\begin{definition} To any parameterized family of (smooth) vector fields $f=f_p(\ubar) \in T_pM$ 
depending upon $\ubar \in \RR$,
one associates the {\bf hyperbolic conservation law on the manifold $(M, \omega)$} 
\be
\label{star}
\del_t u + \dive (f(u))  = 0, \qquad u=u(t,p), \quad t \geq 0, \,  p \in M.  
\ee
The parametrized flux field $f$ is called {\bf geometry-compatible} if constants are trivial solutions, that is 
\be
\label{333}
( \dive f)_p (\ubar) = 0, \qquad \ubar \in \RR, \, p \in M 
\ee
in the distribution sense. A {\bf convex entropy pair} is a convex function  $U: \RR \to \RR$ together with 
a parametrized family of vector field  
$$
F_p(\ubar) := \int^{\ubar} \del_u U(v) \, \del_u f_p(v) \, dv \in T_pM, \qquad \ubar \in \RR, \, p \in M.   
$$ 
\end{definition}

The equation  (\ref{star}) is geometric in nature and does not depend on a particular choice of 
local 
coordinates on $M$.  The flux field is a section of the tangent bundle $TM$ 
and, in general, there is no concept of ``spatially constant flux'', unlike in the Euclidian case
where
one may arbitrarily choose a vector $f_p(\ubar)$ at some point $p \in \RR^n$
and parallel-transport it in $\RR^n$, trivially generating a constant (and, therefore, smooth) vector field. 
The geometry-compatibility condition is natural if one observes that it is 
satisfied by the shallow water equation.

\begin{remark} 
\label{21R} 
An important class of geometry-compatible conservation laws is based 
on {\bf projected gradient flux} fields, as it is called in \cite{BL}, 
and is defined as follows. 
Consider the unit sphere $\Sphere$ embedded in the Euclidian space $\RR^3$ 
and a function $h:\RR^3 \times \RR \to \RR$ being fixed, and 
define  
$$
f_p (\ubar) : = p \wedge (\nabla_{\RR^3} h)_p(\ubar) \in T_p\Sphere, \qquad p \in \Sphere \subset \RR^3.
$$  
\end{remark} 

Now, to formulate the initial value problem associated with (\ref{star}) we are given an initial data $u_0: M \to \RR$ 
and we impose 
\be
\label{starstar}
u(0,p) = u_0(p), \quad p \in M.      
\ee
It is well-known that entropy inequalities must be imposed on weak solutions to conservation laws. 
In this juncture, it is interesting to point out the following property of geometry-compatible flux field: 
If $u=u(t,p)$ is a smooth solution and $U=U(\ubar)$ is an arbitrary function, then  clearly 
$$
\del_t U(u) + \del_u U(u) \, \dive ( f_p(u) ) = 0,  
$$
and it can then be checked that a necessary and sufficient condition for the existence of a vector field $F=F_p(\ubar) \in T_pM$ 
such that 
$$
\del_u U(v) \, \dive ( f_p(v) ) = \dive ( F_p(v) ) 
$$
for all functions $v=v(p)$ is precisely that the flux $f_p(\ubar)$ be geometry-compatible. This observation motivates the following definition which impose the entropy inequality 
\be
\label{294} 
\del_t U(u) + \dive ( F(u) ) \leq 0
\ee
 in the sense of distributions for all convex $U$.  
  
\begin{definition} An {\bf entropy solution to the conservation law} (\ref{star}), 
posed on the manifold $(M,\omega)$ and assuming the initial data $u_0 \in L^\infty(M)$, 
is a scalar field $u \in L^\infty\bigl(\RR^+ \times M\bigr)$ satisfying 
the entropy inequalities 
$$
\int_{\RR^+}\int_M \bigl( U(u) \, \del_t \theta
             + d\theta ( F(u) ) \bigr) \, \omega dt 
 + \int_M  U(u_0) \, \theta(0,\cdot) \, \omega \geq 0 
$$
for all convex $(U,F)$ and test-functions $\theta \geq 0$.
\end{definition}

Note that integration is performed w.r.t.~the volume form on $M$, and  
the differential $d\theta$ is a $1$-form field acting on vector fields.
Our main existence result is as follows. 
 
\begin{theorem}[Well-posedness theory on a manifold] 
Let $(M,\omega)$ be a compact $n$-manifold endowed with an $L^\infty$ volume form, and 
$f=f_p(\ubar) \in T_pM$ be a geometry-compatible flux field.  

\begin{enumerate}

\item Given any data $u_0 \in L^\infty(M)$, the initial value problem associated with the conservation law
(\ref{star}) 
admits a unique entropy solution $u \in L^\infty(\RR^+ \times M)$. 
  
\item Entropy solutions satisfy the following $L^q$ stability property for all $q\in [1, \infty]$ and $t_2 \geq t_1$ 
\be
\label{498}
 \|u(t_2) \|_{L_\omega^q(M)} \leq \|u(t_1) \|_{L^q_\omega(M)}. 
\ee

\item Moreover, the following $L^1$ semi-group property holds for any two entropy solutions $u,v$ and $t_2 \geq t_1$ 
\be
\label{499}
\| v(t_2) - u(t_2) \|_{L^1_\omega(M)} \leq \| v(t_1) - u(t_1) \|_{L_\omega^1(M)}. 
\ee
\end{enumerate} 
\end{theorem}

This theorem provides a generalization of Kruzkov's theory to a manifold, endowed with a volume form with 
possibly low regularity. The $L^q$ stability and 
$L^1$ contraction properties provide {\sl geometry-independent} bounds, which are useful to design 
approximation schemes.  For further discussions on the well-posedness theory based on vector fields, 
we refer the reader to Ben-Artzi and LeFloch~\cite{BL} and
 Amorim, Ben-Artzi and LeFloch~\cite{ABL}, 
as well as Panov~\cite{Panov1,Panov2}.

Consider a geometry-compatible conservation law with flux field $f$ 
together with some initial data $u_0 \in L^\infty(M)$.
Then, an {\bf entropy measure-valued solution} to (\ref{star})-(\ref{starstar}), by definition,
 is a measure-valued field 
$(t,p) \in \RR^+\times M \mapsto \nu_{t,p} \in \mbox{Prob}(\RR)$ such that 
$$
\int_{\RR^+} \int_M  \Bigl( \la \nu_{t,p}, U \ra \, \del_t \theta 
+ \la \la \nu_{t,p}, F \ra, d\theta \ra \Bigr) \, \omega(p) dt 
+ \int_M  U(u_0) \, \theta(0,\cdot) \, \omega \geq 0 
$$
for all convex entropy pairs and test-functions $\theta= \theta(t,p) \geq 0$. 
Following DiPerna \cite{DiPerna} we can establish the following uniqueness result. 

\begin{theorem}[Uniqueness and compactness framework]  
If $\nu$ is an entropy measure-valued solution to the conservation law 
(\ref{star}) posed on the manifold $(M,\omega)$ and 
assuming some initial data $u_0 \in L^\infty(M)$ at $t=0$, 
then 
$$
\nu_{t,p} = \delta_{u(t,p)} \quad  \mbox{ a.e. } p \in M.  
$$
where $u \in L^\infty(\RR^+ \times M)$ denotes the unique entropy solution to the same problem. 
\end{theorem}

This theorem is very useful to establish the convergence of approximation schemes   
since, in comparison with the standard technique based on a total variation bound, 
it requires weaker and {\sl geometrically more natural} bounds on (the entropy dissipation of) 
the approximate solutions, 
as well as {\sl weaker assumptions} on the manifold geometry. 
 
In the Euclidian case, the total variation diminishing (TVD) property was found to be very useful for the development 
of high-order schemes. 
On manifolds, provided $\omega$ is sufficiently smooth and the initial data $u_0$ has bounded variation
then $u(t)$ has bounded variation for all $t \geq 0$. In fact, 
there exists $C>0$ depending on the geometry of $M$ and $\| u_0\|_{L^\infty(M)}$
such that 
$$
TV(u(t)) \leq e^{C \, t} \, ( 1 + TV(u_0) ), 
$$
and the total variation might be unbounded as $t \to \infty$. 

Introduce the total variation along a given vector field $X$ by 
\be
\label{55}
TV_X(u) := \sup_{\| \phi\|_{L^{\infty}}} \int_M u \, \dive ( \phi \, X) \, \omega. 
\ee

\begin{theorem}[Total variation estimate]
Consider a geometry-compatible conservation law (\ref{star}) posed on the manifold $(M,\omega)$ and 
suppose that the initial data have bounded variation. Then, the {\bf total variation diminishing property} (TVD) 
$$
TV_X(u(t_2)) \leq TV_X(u(t_1)), \qquad t_2 \geq t_1, 
$$  
holds provided the flux field $f$ and the given field $X$ satisfy the compatibility condition
\be
\label{56}
\Lcal_X \big( \del_u f(\ubar) \big) := \big[ \del_u f(\ubar), X \big] =0.
\ee
On the other hand, from the time-invariance of the equation and the $L_\omega^1$ contraction property,
it follows that 
\be
\label{666}
\| \dive (f(u(t))) \|_{{\mathcal M}_\omega(M)} := \sup_{\|\theta\|_{C^0} \leq 1} \int_M d\theta(f(u(t))) \, \omega 
\ee
is diminishing in time. 
\end{theorem}

The condition (\ref{56}) is usually not satisfied in the examples of interest, 
except when there exists some decoupling into a one-parameter family of independent, one-dimensional equations. 
(See below.) Hence, the curved geometry restricts the class of flux enjoying the TVD property, and we can not 
expect to rely on a TVD property to develop approximation schemes. 
In contrast, in the Euclidian case no restriction is implied on ($p$-independent) flux and the TVD property
always holds.

\begin{remark} Returning to the projected gradient flux field on $\Sphere \subset \RR^3$ introduced in Remark~\ref{21R}
and choosing a vector field $X$ of the form 
$$ 
X_p := p \wedge (\nabla_{\RR^3} k)_p,  \qquad p \in \RR^3,  
$$
then the condition $\Lcal_X \big( \del_u f (\ubar)\big) = 0$ on $\Sphere$
holds if and only if the vectors $\del_u f_p(\ubar)$ and $X_p$ are parallel, and 
the proportionality factor $C(\ubar, p)$ is constant along the integral curves of $X$: 
$$
\del_u f_p(\ubar) = \del_u C(\ubar, p) \, X_p, 
\quad 
X(C(\ubar, p)) = 0, \qquad \ubar \in \RR, \, p \in M. 
$$ 
\end{remark}


\section{Hyperbolic conservation laws based on differential forms}

In this section based on LeFloch and Okutmustur \cite{LO1,LO2}, we consider 
hyperbolic conservation laws posed on an $(n+1)$-dimensional differentiable manifold $M$, referred to as a spacetime.  
We emphasize that, in the framework now developed, no geometric structure is a priori imposed on $M$. 

\begin{definition}  The {\bf hyperbolic conservation law on the differentiable $(n+1)$-manifold $M$} 
associated with a parametrized family of $L^\infty$ $n$-form fields 
$p \in M \mapsto \omega_p(\ubar)$ depending smoothly upon the parameter $\ubar$, 
by definition, reads 
\be
\label{455}
d( \omega(u)) = 0, 
\ee
whose unknown is the scalar field $u: M \to \RR$. The parametrized flux field $\omega$
 is called {\bf geometry-compatible} if it is closed, that is, for each $\ubar \in \RR$, 
$(d \omega) (\ubar) = 0$  in the distribution sense. 
\end{definition}

By definition, weak solutions $u \in L^\infty(M)$ must satisfy the following weak form of the conservation law 
\be
\label{weak3} 
\int_M d\theta \wedge \omega(u) = 0 \qquad \mbox{ for all test-functions } \theta: M \to \RR. 
\ee 
Then, to pose the initial value problem we impose the following {\bf global hyperbolicity condition:} 
$M$ admits a smooth foliation by compact, oriented hypersurfaces, i.e.,  
$$
M = \bigcup_{t \geq 0} \Hcal_t. 
$$  
Furthermore, denoting by $i_{\Hcal_t}: \Hcal_t \to M$ the injection map, we impose that each hypersurface  
is {\bf spacelike} in the sense that the $n$-forms  
\be
\label{888}
\del_u  \big( i^*_{\Hcal_t} \omega(\ubar) \big) >0   
\ee 
are (positive) $n$-volume forms for each $t \geq 0$ and all relevant $\ubar$.

To formulate the initial value problem associated with (\ref{455}) we are given some initial data $u_0 \in L^\infty(\Hcal_0)$ 
and we search for a scalar field $u\in L^\infty(M)$ satisfying 
\be
\label{data2}
u_{|\Hcal_0} = u_0. 
\ee 
In the present framework, a {\bf convex entropy flux} is now 
a parametrized family of $n$-forms $\Omega = \Omega_p(\ubar)$ such that 
there exists a convex $U:\RR \to \RR$ 
$$
\Omega_p(\ubar) := \int^\ubar \del_u U \,\del_u \omega_p \, du', \qquad p \in M, \, \ubar \in \RR.  
$$
If the flux $\omega$ is geometry-compatible, then smooth solutions satisfy the additional conservation laws 
$$
d ( \Omega(u) ) = 0, 
$$
so that the entropy inequalities take the following form
\be
\label{934}
d ( \Omega(u) ) \leq 0. 
\ee

\begin{definition} 
An {\bf entropy solution to the geometry-compatible conservation law} (\ref{455}) 
is a scalar field $u=u(p) \in L^\infty(M)$ 
satisfying for all convex entropy flux $\Omega = \Omega_p(\ubar)$ and smooth functions $\theta \geq 0$ 
$$  
\int_{M} d\theta \wedge \Omega(u) 
+ 
\int_{\Hcal_0} \theta_{\Hcal_0} \, (i^*\Omega)(u_0) 
 \geq 0.  
$$
\end{definition}

Introduce also Kruzkov's entropy flux fields, defined here as parametrized fields of $n$-forms:  
$$
\Omega(u,v) := \sgn(u-v) \big( \omega(u) - \omega(v) \big).
$$

\begin{theorem}[Well-posedness theory on a spacetime. Geometry-compatible flux] Consider a conservation law with a geometry-compatible flux $\omega= \omega(\ubar)$, 
posed on a globally hyperbolic, spatially compact manifold $M = \bigcup_{t \geq 0} \Hcal_t$.
Then, the initial value problem with data $u_0 \in L^\infty(\Hcal_0)$
admits a unique entropy solution $u \in L^\infty(M)$, 
and for any convex entropy flux $\Omega$ the function 
$t \mapsto \int_{\Hcal_t} i_{\Hcal_t}^* \Omega(u)$ 
is well-defined and non-increasing. 
Moreover, for any two entropy solutions $u,v$
$$
\int_{\Hcal_t} i_{\Hcal_t}^* \Omega(u,v)   
$$
is non-increasing in time. 
\end{theorem}
 
\

\begin{theorem}[Well-posedness theory on a spacetime. General flux] Consider a conservation law with general flux
 field $\omega= \omega(u)$, posed on 
globally hyperbolic, spatially compact $M$ with (past) boundary $\Hcal_0$. 
Consider initial data $u_0$ satisfying $\int_{\Hcal_0} | i_{\Hcal_0} \omega(u_0)| < \infty$. 
Then, there exists a semi-group of entropy solutions $u_0 \mapsto u = Su_0$
and for all smooth functions $\theta \geq 0$ 
\be
\label{999}
\int_{M} d\theta \wedge \Omega(u,v) 
+ 
\int_{\Hcal_0} \theta_{\Hcal_0} \, i_{\Hcal_0}^*\Omega(u_0, v_0) 
 \geq 0.  
\ee
Moreover, for any two spacelike hypersurfaces $\Hcal_2$ in the future of $\Hcal_1$
\be
\label{288}
\int_{\Hcal_2}  i_{\Hcal_2}^* \Omega(u,v) \leq \int_{\Hcal_1} i_{\Hcal_1}^* \Omega(u,v) < \infty.  
\ee
\end{theorem}

\begin{remark} On the one-dimensional torus $M =T^1$,  
the case of general flux fields was covered earlier by LeFloch and N\'ed\'elec~\cite{L1988}
under the assumption that $f$ is strictly convex in $\ubar$. Then, 
entropy solutions satisfy the generalized Lax formula 
$u(t,x) = (\del_u f)^{-1}\Big({x - y(t,x) \over t}\Big)$, in which the function $y=y(t,x)$ 
is determined by a minimization argument over the family of characteristics defined by 
$$
\del_s X = (\del_u f \circ f^{-1}_\pm)\Big( {c \over \omegab(X)} \Big), \qquad X(0) = y, \quad c \in \RR. 
$$
Interestingly, it was established in \cite{L1988} that any two entropy solutions $u,v: \RR^+ \times T^1 \to \RR$ 
satisfy 
$$
\|v(t_2) - u (t_2) \|_{L_\omega^1(T^1)} \leq  \|v(t_1) - u (t_1) \|_{L_\omega^1(T^1)}, \qquad t_2 \geq t_1. 
$$  
\end{remark}


\section{Intrinsic finite volume schemes} 

\subsection{General intrinsic approach}

Using the framework based on $n$-diffe\-rential forms, Amorim, LeFloch, and Okutmustur~\cite{ALO} have introduced 
finite volume schemes, which generate 
piecewise constant approximations on unstructured triangulations. One considers here
the initial value problem associated with the conservation law (\ref{455}) 
and discretize the weak form (\ref{weak3}) of the conservation law. 
Indeed, the finite volume approach applies directly to the
geometric weak form and, in particular, no local coordinates should be chosen at this stage. 
The finite volume discretization on a manifold $M$ requires only the $n$-differential structure
put forward in the previous section. 

One introduces a family of triangulations made of curved polyhedra $K$, each of them having a past boundary $e_K^-$, 
a future boundary $e_K^+$, and ``vertical boundary'' elements $e \in \del K^0$.  For each boundary element $e$ 
there is precisely one element $K_e$ such that $K \cap K_e = e$. 
Then, one defines the notion of ``total flux'' along spacelike hypersurfaces $e_K^\pm$ by setting 
$$
\int_{e_K^\pm} i^* \omega(u) \approx \int_{e_K^\pm} i^* \omega(u_{e_K^\pm}) =: q_{e_K^\pm} (u_{e_K^\pm}),
$$  
where $u: M \to \RR$ denotes the entropy solution of the initial value problem under consideration. 
  
Then, applying Stokes theorem to the conservation law (\ref{455}), we derive the {\bf finite volume scheme} 
\be
\label{485}
q_{e_K^+} (u_{e_K^+}) = q_{e^-_K} (u_{e_K^-}) - \sum_{e\in \del K^0} \qbf_{e,K}(u_{e^-_K},u_{e^-_{K_e}}),
\ee
in which the {\bf total flux} along vertical faces are {\sl scalars}
determined from Lipschitz continuous, numerical flux $\qbf_{e,K} : \RR^2 \to \RR$ satisfying
\begin{enumerate}

 \item Consistency property: 
 $$
 \qbf_{K,e}(\ubar, \ubar )= q_e (\ubar),
$$
 \item Conservation property: 
 $$
 \qbf_{K,e}(\ubar, \vbar) = - \qbf_{K_e, e}(\vbar, \ubar),
$$
 \item Monotonicity property: 
 $$
 \del_{\ubar} \qbf_{K,e}\geq 0, 
       \qquad \del_{\vbar} \qbf_{K,e} \leq 0
       $$
\end{enumerate} 
for all $\ubar, \vbar \in \RR$. For instance, Lax-Friedrich flux or Godunov flux may be used to determine the functions  $\qbf_{e,K}$. 
The use of total flux leads to a geometrically natural setting, in which the volume of 
elements $K$ or boundary elements $e$ do not arise, so that the proposed formulation is both 
conceptually and technically
simpler.

The convergence of finite volume schemes toward the entropy solution of the initial value problem
is established in \cite{ALO,LO2}. Let us mention here the main difficulties:
possible blow-up of the TV norm; derivation of 
uniform estimates tied to the geometry, only involving $n$-dimensional volumes. 
The main issues in our analysis are as follows: 
\begin{enumerate}

\item Derivation of a discrete maximum principle. 

\item Derivation of discrete entropy inequalities. 

\item Derivation of entropy dissipation estimates. 

\item The sequence of finite volume approximation generate an entropy measure-valued solution $\nu_{p}$ 
which, by our generalization of DiPerna's uniqueness theorem, coincides with $\nu_{p} = \delta_{u(p)}$, 
which implies that the scheme converges strongly.  

\end{enumerate}
 
DiPerna's measure-valued solutions were used to establish the convergence of sche\-mes by Szepessy \cite{Szepessy,Szepessy2},
 Coquel and LeFloch \cite{CL1,CL2,CL3}, and Cockburn, Coquel, and LeFloch \cite{CCL00,CCL2}.
For many related results and a review about the convergence techniques for hyperbolic problems,
we refer to Tadmor \cite{Tadmor} and Tadmor, Rascle, and Bagnerini \cite{TRB}.
Further hyperbolic models, including also a coupling with elliptic equations, as well
as many applications were successfully
investigated by Kr\"oners \cite{Kroener}, and Eymard,
Gallouet, and Herbin \cite{EGH}. For higher-order schemes, see the paper by Kr\"oner, Noelle, and Rokyta \cite{KNR}.
An alternative approach to the convergence of finite volume schemes
was discovered by Westdickenberg and Noelle \cite{WN}.


\subsection{A second-order intrinsic scheme on the sphere}  

Ben-Artzi, Falcovitz, and LeFloch \cite{BFL} adopted the fully intrinsic approach and
extended the finite volume scheme to second-order accuracy (by using the generalized Riemann problem methodology),  
while 
implementing it for classes of hyperbolic conservation laws posed on the sphere $\Sphere$.
The sphere is the most important geometry for the applications to geophysical flows. 
Earlier works in this field were not fully geometric, and regarded the manifold as embedded in some Euclidian space
and using projection techniques from the ambient space to the sphere.  
A fully intrinsic approach has definite advantages in terms of convergence and accuracy, and 
at the first-order at least, rigorous convergence analysis can be made. 

We present here some features of the proposed scheme and refer to \cite{BFL} for further discussions, especially
on the implementation.  We consider the conservation law 
\be
\label{395}
\del_t u(t,x) + \nabla_{\Sphere} \cdot (F(x,u(t,x))) = 0
\ee
with flux vector tangent to the sphere, in the general form
$$
F(x,u) = n(x) \wedge \Phi(x,u). 
$$
We are particularly interested in geometry-compatible flux vectors, and have 
proposed the following broad class of {\bf gradient flux vector} 
$$
\Phi(x,u) = \nabla h(x,u). 
$$
In particular, this includes  {\bf homogeneous flux vectors} having $\Phi=\Phi(u)$. 
Interestingly enough, even in this case, the corresponding flux covers a broad and non-trivial class of interest. 

Several classes of solutions of particular interest have been constructed, which exhibit 
very rich wave structure:

\begin{enumerate}

\item Spatially periodic solutions.

\item Non-trivial steady state solutions.

\item Confined solutions that remain compactly supported for all times.  

\end{enumerate}

Explicit formulas for such solutions together with various choices of corresponding flux vector fields 
 can be found in \cite{BFL}. These formulas are useful for numerical investigations.

An intrinsic triangulation was used, a {\sl web-like mesh} made of segments of longitude and latitude lines. 
Artificial coordinate-singularities at the poles are coped with by suitably adapting the mesh in a non-conformal fashion as one approaches the North and South poles. 
To derive explicit formulas for the flux, for instance, in the case of an homogeneous flux vector $\Phi(u)$, one may 
introduce the following basis of the tangent space
$$
i_\lambda = -\sin\lambda         \,i_1
                 +\cos\lambda         \,i_2,\
                 \qquad 
 i_\phi  = 
                 -\sin\phi\,\cos\lambda\,i_1
                 -\sin\phi\,\sin\lambda\,i_2
                 +\cos\phi            \,i_3; 
$$
in terms of some Euclidian basis $i_1, i_2, i_3 ±in \RR^3$. Then, by setting 
$$
\Phi(u) = f_1(u)\, i_1 + f_2(u)\, i_2 + f_3(u)\, i_3,
$$
one finds 
$$
F(x,u) = F_\lambda(\lambda,\phi,u) \, i_\lambda + F_\phi(\lambda,\phi,u) \, i_\phi, 
$$
with
\be
\label{493}
F_\lambda(\lambda,\phi,u) =  
                               f_1(u)\,\sin{\phi}\,\cos{\lambda}
                              +f_2(u)\,\sin{\phi}\,\sin{\lambda}
                              -f_3(u)\,\cos{\phi},
\ee
and 
\be
\label{947}
F_\phi(\lambda,\phi,u) = -f_1(u)\,\sin{\lambda}
                 +f_2(u)\,           \cos{\lambda}.
\ee

The intrinsic finite volume discretization of the conservation law is then based on Stokes formula applied 
$$
\dive_{S^2} F =  
{1 \over \cos \phi} \Big( {\del \over \del \lambda} F_\lambda + {\del \over \del \phi} (F_\phi \, \cos \phi ) \Big). 
$$
The flux are computed by explicit integration along longitude curves, or latitude curves. Length and areas are 
computed exactly, which allows one to derive a discrete version of the geometric compatibility condition.


\section{Perspectives and open problems} 

In conclusion, based on our theoretical investigations (reported in this brief review),  
a shock-capturing numerical method was designed and implemented in \cite{BFL}, i.e.~a 
fully intrinsic version of the second-order Godunov-type, 
finite volume scheme.  
This version, in its principle at first-order accuracy at least, extends to systems 
like the shallow water equations on the sphere with topography. 
Furthermore, the associated discrete properties (conservation, geometric compatibility, asymptotic-preserving) 
established in \cite{BFL} extend to systems as well. 

In the setup advocated here for the design numerical schemes, 
geometric terms in hyperbolic conservation laws are
taken into account in a direct way by discretizing the geometric (weak) formulation of the conservation law, 
instead of using local coordinates and regarding the hyperbolic equation as 
an equation with variables coefficients and source-terms,
or instead of viewing the manifold
as embedded in a higher dimensional Euclidian space.  
This fully intrinsic approach was validated, both theoretically and numerically. 
It will be interesting to further investigate the properties of these finite volume schemes 
on general curved geometries, especially in the context of complex flows and geometries and 
to study the large-time asymptotics of solutions. 

In the shallow water equations on the sphere, one unknown ---the mass density--- is a scalar field and can be treated as 
explained in the present work. However, the other unknown ---the velocity--- is a vector field, and it remains
a challenging open problem to carry out the general Riemann problem methodology in this vector- or, more generally,
tensor-valued context. This is an important issue,
since understanding large-scale atmospheric and oceanic motions depends upon designing robust and accurate numerical 
techniques that, we advocate, should directly incorporate curved geometric effects.

The class of geometry-compatible conservation laws on spacetimes forms  
an analogue of the inviscid Burgers equation, which has served as an important simplified model 
for the development of shock-capturing methods.  
The treatment of boundary conditions on non-compact manifolds is an important issue, tackled 
in \cite{GL}, while error estimates a la Kuznetsov were derived in \cite{ALN,LNO}.


\bibliographystyle{plain}

\end{document}